\documentclass[a4paper,10pt]{article}
\usepackage[utf8]{inputenc}
\usepackage{amssymb,latexsym}
\usepackage{amsmath,amsthm}
\usepackage{graphicx}
\usepackage{color,epsfig}
 
\newtheorem{theorem}{Theorem}
\newtheorem{lemma}{Lemma}

\def\qed{\ifhmode\unskip\nobreak\fi\quad\ifmmode\Box\else$\Box$\fi}

\def\cro{{\mbox {\sc cr}}}

\def\pcr{{\mbox {\sc pcr}}}
\def\blue{{\mbox {\sc blue}}}
\def\red{{\mbox {\sc red}}}

\title{A slightly better bound on the crossing number in terms of the pair-crossing number}
\author{J\'anos Karl\\
\small  Budapest University of Technology and Economics \\[-0.8ex]
\small \texttt{karlj@math.bme.hu}\\
and\\
G\'eza T\'oth\thanks{Supported by National Research, Development and Innovation Office, NKFIH,
K-131529 and ERC Advanced Grant "GeoScape" 882971.}\\
\small Alfr\'ed R\'enyi Institute of Mathematics, Budapest
\small and
\small Budapest University of Technology and Economics \\[-0.8ex]
\small \texttt{geza@renyi.hu}}

\begin{document}

\maketitle

\begin{abstract}
The {\em crossing number} of a graph $G$, $\cro(G)$, is the minimum number of crossings,  
the {\em pair-crossing number}, $\pcr(G)$, is the minimum number of pairs of crossing edges
over all drawings of $G$.  
In this note we show that 
$\cro(G)=O(\pcr(G)^{3/2}\log\pcr(G))$, which is an improvement of the result of Matou\v sek, by a log factor.
\end{abstract}

\section{Introduction}

By a {\em graph} we always mean a {\em simple graph}, that is, a graph with no loops and parallel edges.
We use the term {\em multigraph} if loops and parallel edges are allowed.
A {\em drawing} of a (multi)graph in the plane is a representation
such that vertices are represented by  
distinct points and its edges by curves connecting 
the corresponding points. We 
assume that no edge passes through any vertex
other than its endpoints, no two edges touch each other
(i.e., if two edges have a common interior point, then at this
point they properly cross each other), no three edges
cross at the same point, and two edges cross only finitely many times. 

The {\em crossing number} of a graph $G$, $\cro(G)$, is the minimum number of crossings 
(crossing points)
over all drawings of $G$. 
The {\em pair-crossing number}, $\pcr(G)$, is the minimum number of pairs of crossing edges
over all drawings of $G$.  
In an optimal drawing for $\cro(G)$, any two edges cross at most once. Therefore, 
it is not easy to see the difference between these two definitions.
Indeed, there was some confusion in the literature between these two notions, 
until the systematic study of their relationship \cite{PT00}, \cite{S17}. 
Clearly, $\pcr(G)\le\cro(G)$, and 
in fact, we cannot rule out the possibility, that $\cro(G)=\pcr(G)$ 
for every graph $G$. Probably 
it is the most interesting open problem in this area.
From the other direction, the best known bound is 
$\cro(G)=O(\pcr(G)^{3/2}\log^2\pcr(G))$ \cite{M14}.
In this note we slightly improve it.

\begin{theorem}\label{pcr-cr}
  For any graph $G$,
  $\cro(G)=O(\pcr(G)^{3/2}\log\pcr(G))$.
\end{theorem}

\section{Proof of Theorem \ref{pcr-cr}}

A string graph is the intersection graph of continuous arcs in the plane. 
Vertices of the graph correspond to continous curves (strings)
in the plane such that
two vertices are connected by an edge if and only if the corresponding strings
intersect each other.

Suppose that $G(V, E)$ is a graph of $n$ vertices.
A {\em separator} in a graph $G$ is subset $S\subset V$
for which there is a partition $V=S\cup A\cup B$, $|A|, |B|\le 2n/3$, and there
is
no edge between $A$ and $B$.
According to the Lipton-Tarjan separator theorem, \cite{LT79}, every planar
graph
has a separator of size $O(\sqrt{n})$. 
This result has been generalized in several directions,
in particular, for string graphs, by Fox and Pach  \cite{FP10}.


\begin{theorem}\label{FP10separetor}
{\rm \cite{FP10}} 
There is a constant $c$ such that 
for any string graph $G$ with $m$ edges, 
there is a separator of size at most $cm^{3/4}\sqrt{\log m}$.
\end{theorem}

The bound  has been improved by Matou\v sek \cite{M14} to  $cm^{1/2}{\log m}$ and then
by Lee \cite{L16} to  $cm^{1/2}$, which is asymptotically optimal.

\begin{theorem}\label{L16separator}
{\rm \cite{L16}} 
There is a constant $c$ such that 
for any string graph $G$ with $m$ edges, 
there is a separator of size at most $cm^{1/2}$.
\end{theorem}

Using the separator theorem of Fox and Pach,
T\'oth \cite{T13} proved that for any graph $G$,
$$\cro(G)=O(\pcr(G)^{7/4}\log^{3/2}(\pcr(G))).$$
Matou\v sek \cite{M14} used the same approach, and  his stronger separator theorem, and proved that
$$\cro(G)=O(\pcr(G)^{3/2}\log(\pcr(G))).$$
To prove Theorem \ref{pcr-cr}, we simply apply the separator theorem of
Lee \cite{L16} with the same method.

\medskip

\noindent {\bf Proof of Theorem \ref{pcr-cr}}.
Let $c$ be the constant from the theorem of Lee \cite{L16}.
In a drawing $\cal D$
of a graph $G$ in the plane, 
call those edges which participate in a crossing 
{\em crossing edges},
and those which do not participate in a crossing 
{\em empty edges}.

\begin{lemma}\label{atrajzolas}
Suppose that $\cal D$ is a drawing of a graph $G$ 
in the plane with 
$l>0$ crossing edges and $k\ge 0$ crossing pairs of edges.
Then $G$ can be redrawn such that (i) empty edges are drawn the same way as before, 
(ii) crossing edges are drawn in the neighborhood of the original crossing edges, and
(iii) there are at most $4ck^{3/2}\log l$ edge crossings.
\end{lemma}
  
\medskip

\noindent {\bf Proof of Lemma.}
The proof is by induction on $l$. 
For $l= 0$ the statement is trivial.
Suppose that the statement has been proved for 
all pairs $(l', k')$, where $l'<l$ and consider a drawing of $G$
with $k$ crossing pairs of edges, such that $l$ edges participate in a
crossing. Obviously, ${l\choose 2}\ge k$, and $2k\ge l$, therefore, $2k\ge l>\sqrt{k}$.

We define a string graph $H$ as follows. The vertex set $\overline{F}$ 
of $H$ corresponds to 
the crossing edges of $G$. Two vertices are connected by an edge if 
the corresponding edges cross each other. Note that the endpoints do not count;
if two edges do not cross, the correspondig vertices are not connected even if
the edges have a common endpoint.
The graph $H$ is a string graph, it can be represented by the crossing edges of
$G$, as strings, with their endpoints removed. 
It has $l$ vertices, and $k$ edges. By Theorem \ref{L16separator}, 
$H$ has a separator of size $ck^{1/2}$
that is, the vertices can be decomposed into three sets, ${\overline{F}}_0$, 
${\overline{F}}_1$, ${\overline{F}}_2$, such that 
(i) $|{\overline{F}}_0|\le ck^{1/2}$, 
(ii) $|{\overline{F}}_1|, |{\overline{F}}_2|\le 2l/3$,
(iii) there is no edge of $H$ between  ${\overline{F}}_1$ and ${\overline{F}}_2$.

This corresponds to a decomposition of the set of crossing edges
$F$ into three sets, $F_0$, $F_1$, and $F_2$ such that
(i) $|F_0|\le ck^{1/2}$,
(ii) $|F_1|, |F_2|\le 2l/3$,
(iii) in drawing $\cal D$, edges in $F_1$ and in $F_2$ do not cross each
other.

For $i=0, 1, 2$, let $|F_i|=l_i$.
Let $G_1=G(V, E\cup F_1)$ and $G_2=G(V, E\cup F_2)$, then in the drawing $\cal D$
of the graph $G_i$ has $l_i$ crossing edges. Denote by $k_i$ the number of crossing pairs
of edges of $G_i$ in drawing $\cal D$. 
Then we have 
$k_1+k_2\le k$, $l_1, l_2\le 2l/3$,
$l_1+l_2+l_0=l$.

For $i=1, 2$, 
apply the induction hypothesis for $G_i$ and drawing $\cal D$.
We obtain a drawing ${\cal D}_i$ satisfying the conditions of the Lemma:
 (i) empty edges drawn the same way as before, 
(ii) crossing edges are drawn in the neighborhood of the original crossing edges, and
(iii) there are at most $4ck_i^{3/2}\log l_i$ edge crossings.

Consider the following drawing ${\cal D}_3$ of $G$. 
(i) Empty edges are drawn the same way as in $\cal D$, ${\cal D}_1$, and ${\cal D}_2$,
(ii) For $i=1, 2$, edges in $F_i$ are drawn as in ${\cal D}_i$,
(iii) Edges in $F_0$ are drawn as in ${\cal D}$.
Now count the number of edge crossings (crossing points) in the drawing ${\cal D}_3$.
Edges in $E$ are empty, 
edges in $F_1$ and in $F_2$ do not cross each other,
there are at most 
$4ck_i^{3/2}\log l_i$ crossings among edges in $F_i$. 
The only problem is that edges in $F_0$ might cross edges in $F_1\cup F_2$ and
each other 
several times, so we can not give a reasonable upper bound for the number of
crossings of this type.
Color edges in $F_1$ and $F_2$ blue, edges in $F_0$ red.
For any piece $p$ of an edge of $G$, let 
 $\blue(p)$ (resp. $\red(p)$)
denote  the number of crossings on $p$ with {\em blue}
(resp. {\em red}) edges of $G$.
We will apply the following transformations.

%
%

\smallskip

{\sc ReduceCrossings$(e, f)$}
Suppose that two crossing edges, $e$ and $f$ cross twice, say,
in $X$ and $Y$.
Let $e'$ (resp. $f'$) be the piece of $e$ (resp. $f$) between $X$ and $Y$.
If $\blue(e')<\blue(f')$, or $\blue(e')=\blue(f')$ and 
$\red(e')\le \red(f')$, then redraw $f'$
along $e'$ from $X$ to $Y$. 
Otherwise, redraw $e'$
along $f'$ from $X$ to $Y$. 
See Figure 3.

\smallskip

\begin{figure}[ht]
\begin{center}
\scalebox{0.6}{\includegraphics{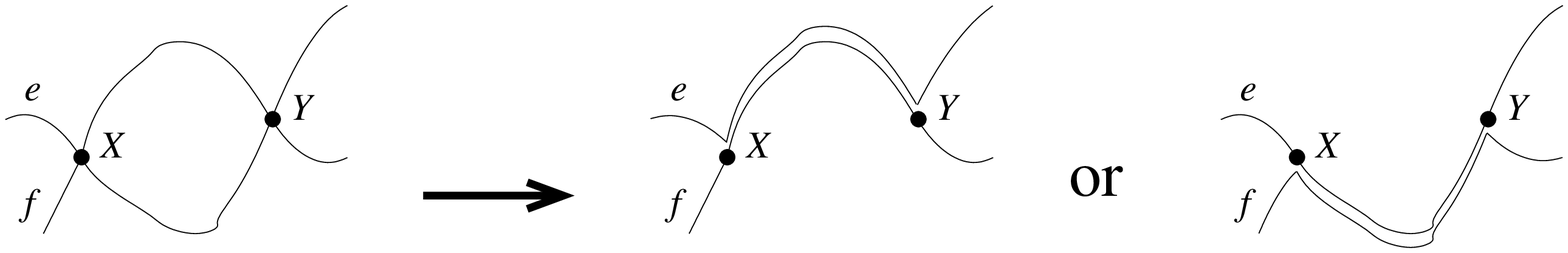}}

\bigskip

\bigskip

\scalebox{0.6}{\includegraphics{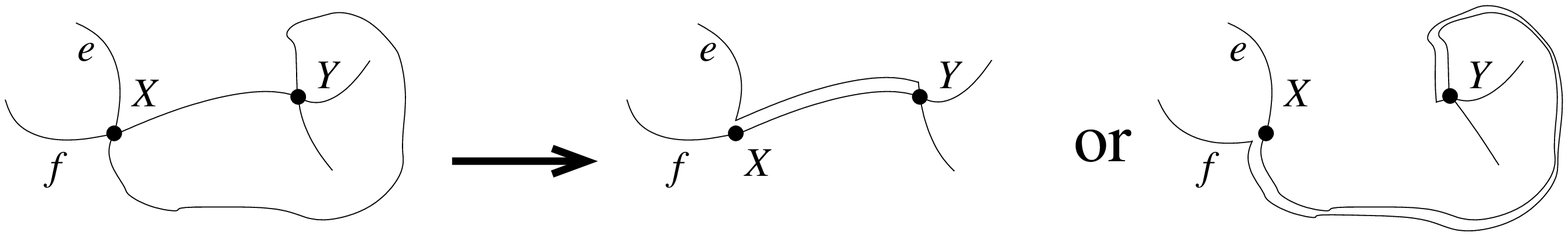}}
\caption{{\sc ReduceCrossings$(e, f)$}}
\end{center}
\end{figure}

\smallskip

Observe that {\sc ReduceCrossings} 
might create self-crossing edges, so we
need another transformation. 

\smallskip

{\sc RemoveSelfCrossings$(e)$}
Suppose that an edge $e$ crosses itself in $X$. Then $X$ appears
twice on $e$. Remove the part of $e$ between the first and last appearance 
of $X$.

\smallskip

Start with drawing ${\cal D}_3$ of $G$, and 
apply 
{\sc ReduceCrossings} and {\sc RemoveSelfCrossings} recursively, 
as long as there are two crossing edges
that cross at least twice, or there is a self-crossing edge. 

Let $BB$, (resp. $BR$, $RR$) denote the number 
of blue-blue (resp. blue-red, red-red) crossings in the current drawing of
$G$. 
Observe, that the triple $(BB, BR, RR)$ lexicographically decreases
with each of the transformations. Indeed, 

\begin{itemize}

\item if $e$ and $f$ are both blue edges then  
{\sc ReduceCrossings$(e, f)$} decreases $BB$,

\item if $e$ is blue and $f$ is red then  either $BB$ decreases,
or if it stays the same then $BR$ decreases,

\item if $e$ and $f$ are both red edges then  
$BB$ stays the same, and either $BR$ decreases,
or if it also stays the same then $RR$ decreases,

\item if $e$ is blue then {\sc RemoveSelfCrossings$(e)$}
decreases $BB$, 

\item and finally, if $e$ is red then $BB$ does not change, $BR$ does not
increase, and $RR$ decreases.

\end{itemize}

Therefore, after finitely many steps we arrive to a drawing 
${\cal D}_4$ of $G$, where any two edges cross at most once, 
and $(BB, BR, RR)$ is lexicograhically not larger than originally.
That is, in the drawing ${\cal D}_4$, 
$BB\le 4ck_1^{3/2}\log l_1 + 2ck_2^{3/2}\log l_2$, 
and any two edges cross at most once, therefore, 
$BR+RR\le l_0l$. 
So, for the total number of crossings we have
$$ 4ck_1^{3/2}\log l_1 + 4ck_2^{3/2}\log l_2+l_0l$$
$$\le 
4ck_1^{3/2}\log (2l/3) + 4ck_2^{3/2}\log(2l/3)+l_0l$$
$$= 
4c(k_1^{3/2}+k_2^{3/2})(\log l+\log (2/3)) + l_0l$$
$$\le 4ck^{3/2}(\log l+\log (2/3)) + l_0l$$
$$\le 4ck^{3/2}\log l - 2ck^{3/2} + l_0l$$
$$\le 4ck^{3/2}\log l - 2ck^{3/2}  + 2ck^{3/2}$$
$$=4ck^{3/2}\log l.$$
$\Box$

Now consider a graph $G$ and let $\pcr(G)=k$. Take a drawing of $G$ with
exactly $k$ crossing pairs of edges. Let $l$ be the total number of crossing
edges.
By Lemma \ref{atrajzolas}, $G$ can be redrawn with at most $4ck^{3/2}\log l$ 
crossings.
Since $2k\ge l$, 
$\cro(G)\le 4ck^{3/2}\log l< 8ck^{3/2}\log k$. This concludes the
proof of Theorem \ref{pcr-cr}. $\Box$

\end{document}